\documentclass[
11pt,%
final,%
tightenlines,%
notitlepage,%
twoside,%
onecolumn,%
nobibnotes,%
nofootinbib,%
superscriptaddress,%
noshowpacs,%
noshowydk,%
showdoi,%
centertags,%
maik]%
{revtex4n}

\input jpack_e
\setcaptionmargin{5mm}
\onelinecaptionstrue

\columntovsizetrue
\allowdisplaybreaks

\begin{document}

\newtheorem{teo}{Theorem}
\newtheorem*{teon}{Theorem}
\newtheorem{lem}{Lemma}
\newtheorem*{lemn}{Lemma}
\newtheorem{prp}{Proposition}
\newtheorem*{prpn}{Proposition}
\newtheorem{ass}{Assertion}
\newtheorem*{assn}{Assertion}
\newtheorem{assum}{Assumption}
\newtheorem*{assumn}{Assumption}
\newtheorem{stat}{Statement}
\newtheorem*{statn}{Statement}
\newtheorem{cor}{Corollary}
\newtheorem*{corn}{Corollary}
\newtheorem{hyp}{Hypothesis}
\newtheorem*{hypn}{Hypothesis}
\newtheorem{con}{Conjecture}
\newtheorem*{conn}{Conjecture}
\newtheorem{dfn}{Definition}
\newtheorem*{dfnn}{Definition}
\newtheorem{problem}{Problem}
\newtheorem*{problemn}{Problem}
\newtheorem{notat}{Notation}
\newtheorem*{notatn}{Notation}
\newtheorem{quest}{Question}
\newtheorem*{questn}{Question}

\theorembodyfont{\rm}
\newtheorem{rem}{Remark}
\newtheorem*{remn}{Remark}
\newtheorem{exa}{Example}
\newtheorem*{exan}{Example}
\newtheorem{cas}{Case}
\newtheorem*{casn}{Case}
\newtheorem{claim}{Claim}
\newtheorem*{claimn}{Claim}
\newtheorem{com}{Comment}
\newtheorem*{comn}{Comment}

\theoremheaderfont{\it}
\theorembodyfont{\rm}

\newtheorem{proof}{Proof}
\newtheorem*{proofn}{Proof}

\selectlanguage{english}
\Rubrika{\relax}
\CRubrika{\relax}
\SubRubrika{\relax}
\CSubRubrika{\relax}
%

\def\JournalNumber{0}
\def\JournalVolume{00}
%
%
%
\nameVolumeRus{}
\CnameVolumeRus{}
\nameIssueRus{\No}
\CnameIssueRus{}
\namePartRus{}
\namePagesRus{}
\nameYearShortRus{}
\JournalNameRus{}
\TranslitJournalNameRus{}
\JournalName{Regular and Chaotic Dynamics}
\JournalISSNCode{1560-3547}
\IssuePrice{}
\TransYearOfIssue{0000}
\TransCopyrightYear{2016}%
\OrigYearOfIssue{}
\OrigCopyrightYear{2016}%
\OrigIssueNo{\JournalNumber}
\OrigVolumeNo{\JournalVolume}
\TransVolumeNo{\JournalVolume}
\TransIssueNo{\JournalNumber}
\TransPartNo{}
\SHORTjournalPREFIX{RCD} 
\LONGjournalPREFIX{RegDyn} 
\BatFileName{call make_ps.bat} 
\BatSwitch{3} 
\IssueName{}
\SupplementNumber{}
\PublicationSerialNumberInYear{0}
\PublicationSerialNumberInVolume{0}
\ConditionalIssueDate{"year","month","day","name","type"}
\PagePrefix{}
\JournalISSNonlineCode{}
\JournalISSNCodeRus{}
\JournalISSNonlineCodeRus{}
\VolumeName{}
\IssnoName{none}
\PartnoName{}
\FpageNamepp{}
\FpageNnamep{}
\FpagePrefix{}
\LpageNnamepp{}
\LpageNamep{}
\LpagePrefix{}
\VolumePageNumbering{}
\JournalPubID{}
\FirstJournalPageNumber{}
\LastJournalPageNumber{}
\makeatletter
\def\MAIKlogo{RCD Editorial Office}
\def\maikpraefix{10.0000/S}
\edef\@ContentsHeadLineB{Simultaneous English language translation of the journal is available from \noexpand\MAIKlogo}
\def\Distributed{Distributed worldwide by Springer. }
\def\ArticlePages#1{\relax}
\@ifxundefined\CONT@sw{\@booleantrue\CONT@sw}{}%
\@booleantrue\showPACS@sw%
\@booleantrue\showKEYS@sw %
\@booleantrue\noOrigJournalVersion@sw
\@booleantrue\noOrigVolumeNo@sw
\@booleanfalse\noTransVolumeNo@sw
\makeatother
\input maikdoi %

\beginpaper


\input engnames
\titlerunning{Discrete Lorenz attractors}
\authorrunning{I. I. Ovsyannikov}
\toctitle{On birth of discrete Lorenz  attractors under bifurcations of 3D maps with nontransversal heteroclinic cycles}
\tocauthor{I.\,I.\,Ovsyannikov}
\title{On birth of discrete Lorenz  attractors under bifurcations of 3D maps 
 with nontransversal heteroclinic cycles}
\firstaffiliation{
}%
\articleinenglish 
\PublishedInRussianNo
\author{\firstname{Ivan\,I.}~\surname{Ovsyannikov}}%
\email[E-mail: ]{ivan.i.ovsyannikov@gmail.com; iovsyann@uni-bremen.de}
\affiliation{
University of Bremen, MARUM, Department of Mathematics\\
Bibliothekstrasse 5, 28359 Bremen, Germany}
\affiliation{
Lobachevsky State University of Nizhny Novgorod, ITMM\\
23 Gagarin av., 603022 Nizhny Novgorod, Russia }%
\begin{abstract}
Lorenz attractors are important objects in the modern theory of chaos. The reason from one side is that they are met in various natural applications (fluid dynamics, mechanics, laser dynamics, etc.). At the same time, Lorenz attractors are robust, in the sense that they are generally not destroyed by small perturbations (autonomous, non-autonomous, stochastic). This allows us to be sure that the observed in the experiment object is exactly the chaotic attractor, rather than a long-time periodic orbit.

Discrete-time analogs of the Lorenz attractor possess even more complicated structure -- they allow homoclinic tangencies of invariant manifolds within the attractor. Thus, discrete Lorenz attractors belong to the class of wild chaotic attractors. These attractors can be born in codimension-three local and certain global (homoclinic and heteroclinic) bifurcations. While various homoclinic bifurcations leading to such attractors were studied, for heteroclinic cycles only cases when at least one of the fixed points is saddle-focus were considered to date. 

In the present paper the case of a heteroclinic cycle consisting of saddle fixed points with a quadratic tangency of invariant manifolds, is considered. It is shown that in order to have a three-dimensional chaos such as the discrete Lorenz attractors, one needs to avoid the existence of lower-dimensional global invariant manifolds. Thus, it is assumed that either the quadratic tangency or the transversal heteroclinic orbit is non-simple. The main result of the paper is the proof that the original system is the limiting point in the space of dynamical systems of a sequence of domains in which the diffeomorphism possesses discrete Lorenz attractors.
\end{abstract}
\keywords{{\em Heteroclinic orbit, rescaling, 3D Henon map, bifurcation, Lorenz attractor}}
\pacs{37C05, 37G25, 37G35}
\received{August 17, 2017}
\revised{Month XX, 20XX}
\accepted{Month XX, 20XX}%
\maketitle

\textmakefnmark{0}{)}%


\section{INTRODUCTION}


Since its discovery, the Lorenz attractor became an important object of the theory of dynamical systems. The reason is that it is a strange attractor that preserves its strangeness under small perturbations. It means that bifurcations may occur inside the attractor, but they do not lead to the appearance of simple attractors e.g. stable periodic orbits. Indeed, in the Lorenz attractor, while changing the parameters, numerous homoclinic bifurcations occur, but they lead only to the appearance and disappearance of saddle periodic orbits.
Based on this, the Lorenz attractor is an object that can be observed in applications and experiments: even when the parameter values and initial conditions are given with  tolerance, one can be sure that the trajectory still converges to a chaotic attractor rather than a periodic orbit with a large period. 

This is not the case for many other known types of attractors, such as the Henon attractor, R\"ossler attractor, attractor in the Chua's circuit. These attractors exist for certain parameter values, but in any neighborhood of these values there exist systems with stable periodic orbits. Such 
attractors are called quasiattractors in the classification by Afraimovich and Shilnikov in \cite{AS83}. The Lorenz attracror belongs to the class of genuine strange attractors according to this classification.


It is known that the Lorenz attractor contains infinitely many saddle periodic orbits, and their stable and unstable manifolds intersect transversely. In paper 
\cite{TS98} an example of a {\em wild spiral attractor} was presented, in which stable and unstable sets in the chaotic attractor can have tangencies. In bifurcations breaking those tangencies, stable periodic orbits still do not appear, so that the attractor is also genuine. The wild spiral attractor belongs to the class of {\em wild hyperbolic attractors}, introduced and described in \cite{TS98}. 

Another representative of this class is the discrete Lorenz attractor (see the definition in \cite{GGOT13}). It can be regarded as a time-discretization of a classical Lorenz attractor under a periodic non-autonomous perturbation \cite{TS08}. They can be met in applications, in particular, in the nonholonomic rattleback model \cite{GGK13} or a two-component convection \cite{EM18}.
Discrete Lorenz attractors are known to be born at local bifurcations
of periodic orbits having three or more
multipliers lying on the unit circle. The following 3D H\'enon map
\begin{equation}
\bar x = y, \;\; \bar y = z, \;\; \bar z = M_1 + B x + M_2 y - z^2
\label{H3D}
\end{equation}
controlled by three independent parameters $M_1$, $M_2$ and $B$ is an example of a model that can undergo such codimension-three bifurcations. 
In  papers \cite{GOST05, GMO06, GGOT13}
it was proved that  map (\ref{H3D}) possesses a discrete Lorenz-like attractor in some open parameter domain near 
point $(M_1=1/4, B=1, M_2 = 1)$, where the map has a fixed point with the triplet $(-1,-1,+1)$ of multipliers. The flow normal form of this bifurcation, after rescaling the coordinates and time, can be brought to the Shimizu-Morioka system:
\begin{equation}  \label{sm}
      \dot X=Y, \;\; \dot Y=X(1-Z)-\lambda Y, \;\;\dot Z=-\alpha Z+X^2.
\end{equation}
This system possesses the Lorenz attractor in some open domain of parameters as proved in \cite{CTZ18}. Then, map (\ref{H3D}) can be regarded as a time-shift map of periodically perturbed system (\ref{sm}). By \cite{TS08}, such a discretization has a chaotic attractor, that was called the discrete Lorenz attractor in \cite{GGOT13}.

This result immediately implies the existence of cascades  discrete Lorenz attractors in global (homoclinic and heteroclinic) bifurcations 
in which map (\ref{H3D}) appears as the first return map. The first such example was considered in \cite{GMO06} 
in case of a quadratic homoclinic tangency to a saddle-focus fixed point with a unit Jacobian. Later analogous results 
were obtained for non-simple homoclinic tangencies to saddle fixed points \cite{GOT14} and homoclinic tangencies to resonant saddle \cite{GO17}. In heteroclinic cycles such bifurcations were studied only when it contains at least one saddle-focus \cite{GST09, GO10, GO13}. Note that the presence of 
saddle-foci in these cases is a very important condition for the existence of Lorenz-like attractors as it prevents 
from the existence of lower-dimensional center manifolds and helps for the dynamics to be effectively three-dimensional (see \cite{T96}).
Another important condition for this is the restriction on the Jacobians in the fixed points. It is based on the fact that
the orbits under consideration may spend unboundedly large time in the neighbourhoods of the saddle fixed point. In the homoclinic case
this means that if the Jacobian differs from one, the phase volumes near such orbits will be either unboundedly expanded or unboundedly contracted,
and the dynamics will have effective dimension less than three. In the same way, for the heteroclinic cases it is necessary to demand for all
the Jacobians not to be simultaneously contracting ($< 1$) or expanding ($> 1$). Thus, in order to get Lorenz attractors in heteroclinic
cycles, one needs to consider ``contracting-expanding'' or ``mixed'' cases.

In the present paper the results of \cite{GST09, GO10, GO13} are extended to the case of heteroclinic cycles that consist of saddles only, i.e. all the fixed points have only real multipliers.
As it is known from \cite{GST93c, Tat01}, in order to have the effective
dimension of the corresponding problem to be equal to $3$,   an additional degeneracy assumptions should be imposed. 
Namely, one of the heteroclinic orbits is {\em non-simple} at the bifurcation moment. The case when the quadratic tangency is non-simple, was also called {\em generalized} tangency in \cite{Tat01}. However, in the heteroclinic cycle there is one more possibility --- when the transversal heteroclinic orbit is non-simple (the required definitions are given below in section~\ref{sec:def}.). Such a case was not studied before.


\begin{rem} The notion of a simple quadratic homoclinic tangency $($that is analogous to the notion of a quasitransversal homoclinic intersection 
{\em \cite{NPT}}$)$ 
was introduced in {\em \cite{GST93c}}. For three-dimensional maps with a homoclinic tangency to a saddle point $O$ with multipliers $\nu_i, i=1,2,3$ such that 
$|\nu_1|<|\nu_2|<|\nu_3|$, it implies the existence of a global two-dimensional invariant manifold $\mathfrak{M}$ for all nearby maps. 
This manifold contains all orbits entirely lying in a small fixed neighbourhood of the homoclinic orbit. In general, it is $C^1$ only and 
particularly hyperbolic. If point $O$ has type $(2,1)$, i.e. $|\nu_{1,2}| < 1 < |\nu_3|$, the manifold is center-stable; if  point $O$ has type $(1,2)$, 
i.e. $|\nu_{1}|<1<|\nu_{2,3}|$, the manifold is center-unstable. It implies that neither periodic nor strange attractors can be born at homoclinic 
bifurcations if $|\nu_2\nu_3|>1$. 
However, as it was shown in {\em \cite{Tat01}}, periodic attractors can appear even in these cases if the homoclinic tangency is non-simple, 
see also  paper {\em \cite{GGT07}} in which the case of the point of type $(1,2)$ was considered in more details. 
These results are very important for the theory of dynamical chaos since they show that the appearance of non-simple homoclinic tangencies can 
destroy the ``strangeness'' of attractors.
\end{rem}

The paper is organised as follows. Section~\ref{sec:def} contains the statement of the problem, main definitions,  three principally
different cases of non-simple heteroclinic orbits are identified, and the main theorem is formulated. In section~\ref{res:het} the first return map is constructed and the rescaling lemma
for all three cases is formulated. Section~\ref{sec:th3proof} contains the proof of the rescaling lemma, from which the statement of the main theorem follows directly.

\begin{figure}[!ht]
\includegraphics[width=8cm]{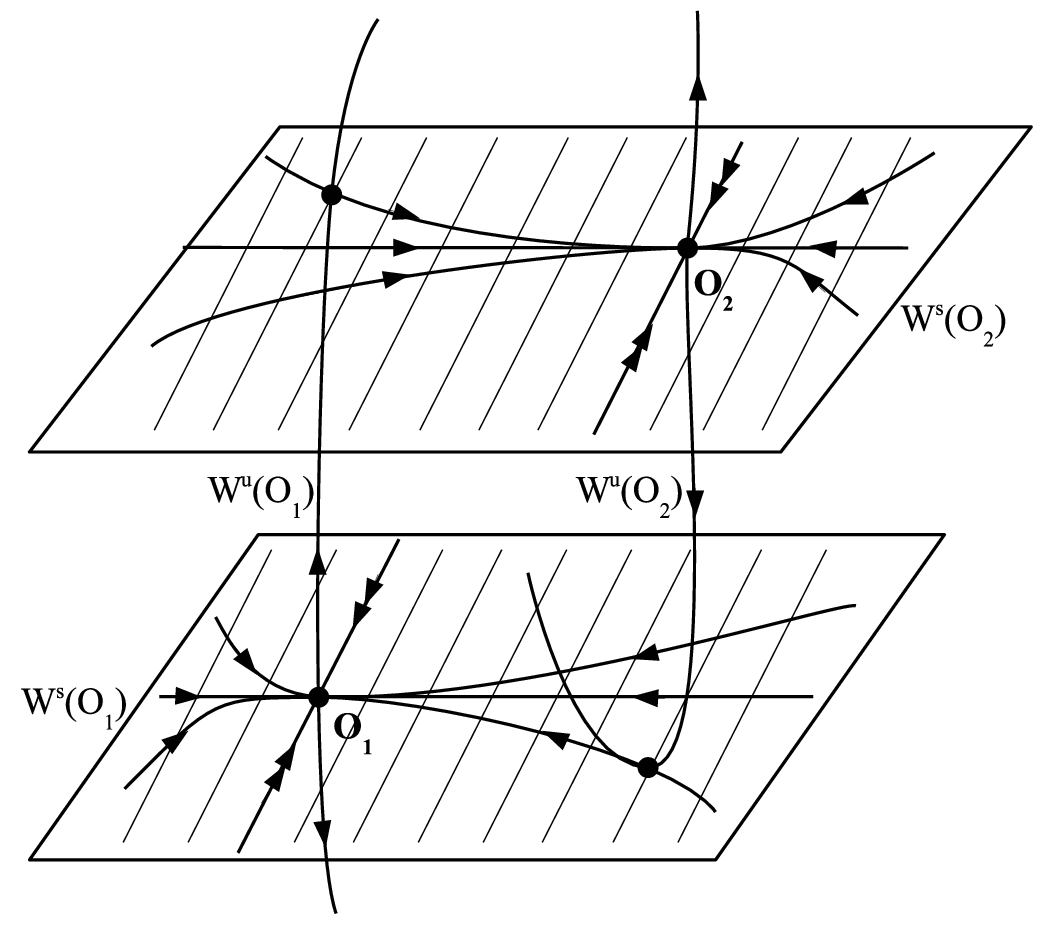} 
\caption{\label{saddles} A heteroclinic cycle consisting of two saddles, with a quadratic tangency of manifolds.} 
\end{figure}



\section{Statement of the problem and main definitions}\label{sec:def}

Consider a three-dimensional orientable $C^r$-diffeomorphism $f_0$, $r \geq 3$, satisfying the following conditions (see Fig.~\ref{saddles}):

{\bf A)} $f_0$ has two fixed points of saddle type: $O_1$ with multipliers $(\lambda_1, \lambda_2, \gamma_1)$ where
    $0 < \lambda_2 < \lambda_1 < 1 < \gamma_1$, and $O_2$ with multipliers $(\nu_1, \nu_2, \gamma_2)$ such that
    $0 < \nu_2 < \nu_1 < 1 < \gamma_2$,

{\bf B)} $J_1 = J(O_1) \equiv \lambda_1 \lambda_2 \gamma_1 < 1$, $J_2 = J(O_2) \equiv \nu_1 \nu_2 \gamma_2 > 1$,

{\bf C)} Invariant manifolds $W^u(O_1)$ and $W^s(O_2)$ intersect transversely at the points of a heteroclinic orbit
$\Gamma_{12}$, the invariant manifolds $W^u(O_2)$ and $W^s(O_1)$ have a quadratic tangency at the points of a heteroclinic orbit
$\Gamma_{21}$.
\\

it is also assumed that the heteroclinic cycle has an additional degeneracy, namely $f_0$ satisfies to one of the following conditions:

{\bf D1)} The transversal intersection of $W^u(O_1)$ and $W^s(O_2)$ is {\em simple}\footnote{The detailed definition of simple and non-simple
heteroclinic orbits will be given below in this section.} and the quadratic 
tangency of $W^u(O_2)$ and $W^s(O_1)$ is {\em non-simple}.

{\bf D2)} The transversal intersection of $W^u(O_1)$ and $W^s(O_2)$ is {\em non-simple} and the quadratic 
tangency of $W^u(O_2)$ and $W^s(O_1)$ is {\em simple}.

The goal of the paper is the study of bifurcations of single-round periodic orbits in generic unfoldings of $f_0$. For this purpose the necessary number of parameters to take, should be identified.
Diffeomorphisms close to $f_0$ and satisfying either conditions {\bf A}--{\bf D1} or {\bf A}--{\bf D2} compose 
locally connected bifurcation surfaces of codimension two in the space of $C^r$-diffeomorphisms, 
thus the number of parameters must be at least two. It is natural to choose 
the splitting parameter of invariant manifolds $W^u(O_2)$ and $W^s(O_1)$ with respect to some 
point of $\Gamma_{21}$ as the first parameter $\mu_1$. 
The second parameter $\mu_2$ is taken to control conditions {\bf D1} or {\bf D2} in such a way that for $\mu_1 = 0$ and $\mu_2 \neq 0$
the corresponding degeneracy disappears i.e. 
in the case {\bf D1} the tangency of $W^u(O_2)$ and $W^s(O_1)$ becomes simple and
in the case {\bf D2} the  intersection of $W^u(O_1)$ and $W^s(O_2)$ becomes simple.
Also note that due to condition {\bf B} we have the contracting-expanding (or mixed) case, which requires one more 
parameter $\mu_3$ that will control the values of the Jacobians 
$J_1$ and $J_2$. It is well known \cite{GST09, GO10, GO13} that the following value effectively plays this role:
\begin{equation}
\label{eq:e3}
\mu_3 = S(f_\mu) - S(f_0), 
\end{equation}
where $S(f)$ is a functional defined as $\displaystyle S(f) = - \frac{\ln J_1}{\ln J_2}$.

In order to define simple and non-simple heteroclinic orbits,
recall some facts from the normal hyperbolicity theory.
Let $O$ be a saddle fixed points of type $(2, 1)$ and $U_0$ be some small neighbourhood of it.
It is known \cite{GST07, GS90, GS92, book} that
diffeomorphism $\left. f_\mu \right|_{U_0}$ for each small $\mu$ can be represented in some $C^r$-smooth local
coordinates $(x_1, x_2, y)$ as follows (the so-called {\em main normal form}):
\begin{equation}
\begin{array}{l}
\bar x_1 \; = \; \lambda_1(\mu) x_1 +
\tilde H_2(y, \mu)x_2 + O(\|x\|^2|y|)  \\
\bar x_2 \; = \; \lambda_2(\mu) x_2 + \tilde R_2(x, \mu) +
\tilde H_4(y, \mu)x_2 + O(\|x\|^2|y|) \\
\bar y \; = \; \gamma(\mu) y + O(\|x\||y|^2) ,\\
\end{array}
\label{t0norm}
\end{equation}
where $\tilde H_{2,4}(0,\mu) = 0$, 
$\tilde R_{2}(x, \mu) = O(\|x\|^2)$.
In coordinates (\ref{t0norm}) the invariant manifolds of saddle fixed point $O$ are 
locally straightened: stable $W^s_{loc}(O): \{ y = 0\}$,
unstable $W^u_{loc}(O): \{ x_1 = 0, \; x_2 = 0\}$ and strong stable $W^{ss}_{loc}(O):  \{ x_1 = 0, \; y = 0\}$.

\begin{figure}[!ht]
\includegraphics[height=6cm]{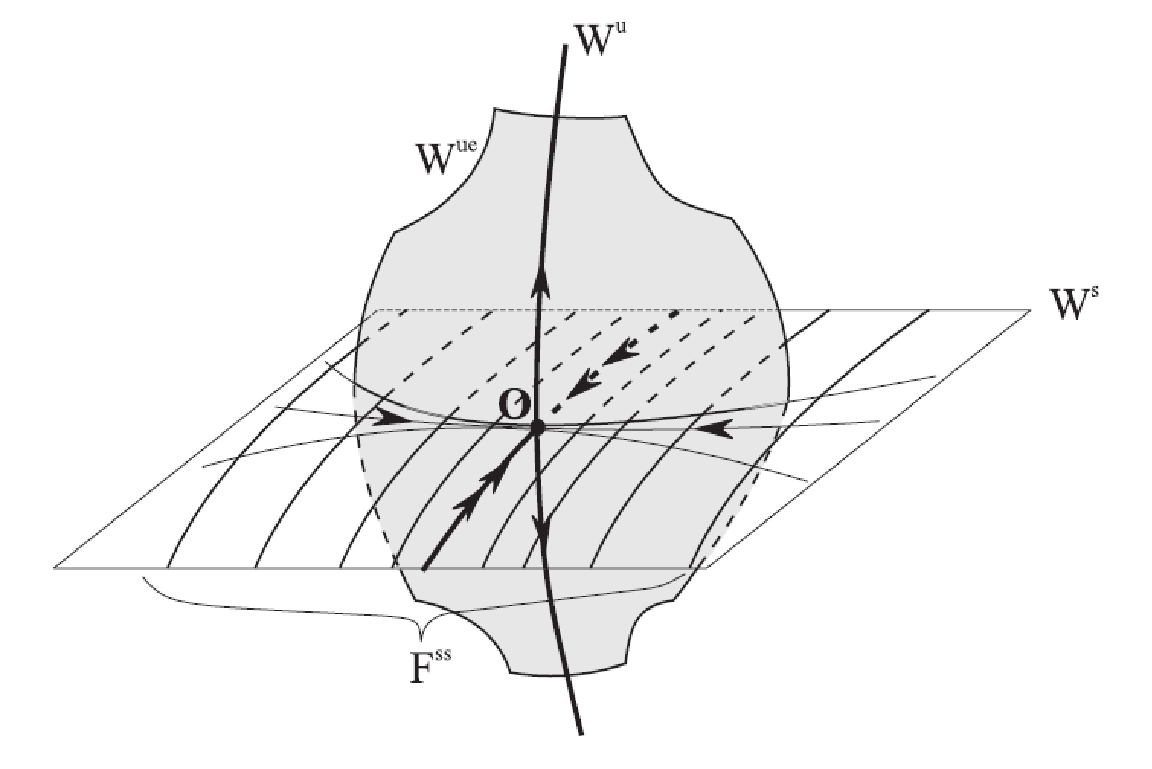}
 \caption{\label{fig04} Invariant structures near a saddle fixed point. A part of the strong stable
foliation $F^{ss}$ containing the strong stable manifold $W^{ss}$
and a piece of one of the extended unstable manifolds $W^{ue}$ containing
$W^u$ and being transversal to $W^{ss}$ at $O$.} 
\end{figure}

According to \cite{book, HPS}, an important role in dynamics is played by an {\em extended unstable manifold }  $W^{ue}(O)$,
see Fig.~\ref{fig04}. By definition, it is a two-dimensional invariant manifold, tangent to the leading stable direction 
(corresponding to $\lambda_1$) at the saddle point and containing unstable manifold $W^u(O)$. Unlike the previous
ones, the extended unstable manifold is not uniquely defined and its smoothness is, generally speaking, 
only $C^{1 + \varepsilon}$. Locally, $W^{ue}_{loc}(O) = W^{ue}(O) \cap U_0$, and the
equation  of $W^{ue}_{loc}(O)$ has the form $x_2 = \varphi(x_1, y)$, where $\varphi(0, y) \equiv 0$ and $\varphi'_{x_1}(0, 0) = 0$.
Note that despite the fact that $W^{ue}(O)$ is non-unique, all of them have the same tangent plane at each point of $W^u(O)$.

Another essential fact is the existence of the
{\em strong stable invariant foliation}, see Fig.~\ref{fig04}.
In $W^s(O)$ there exists a one-dimensional strong stable invariant submanifold $W^{ss}(O)$, which is
$C^r$--smooth and touches at $O$ the eigenvector corresponding to the strong
stable (nonleading) multiplier $\lambda_2$. Moreover, manifold $W^s(O)$ is foliated near $O$ by
the leaves of invariant foliation $F^{ss}$ which is $C^r$-smooth, unique and contains
$W^{ss}(O)$ as a leaf.

\begin{figure}[!ht]
\includegraphics[height=7cm]{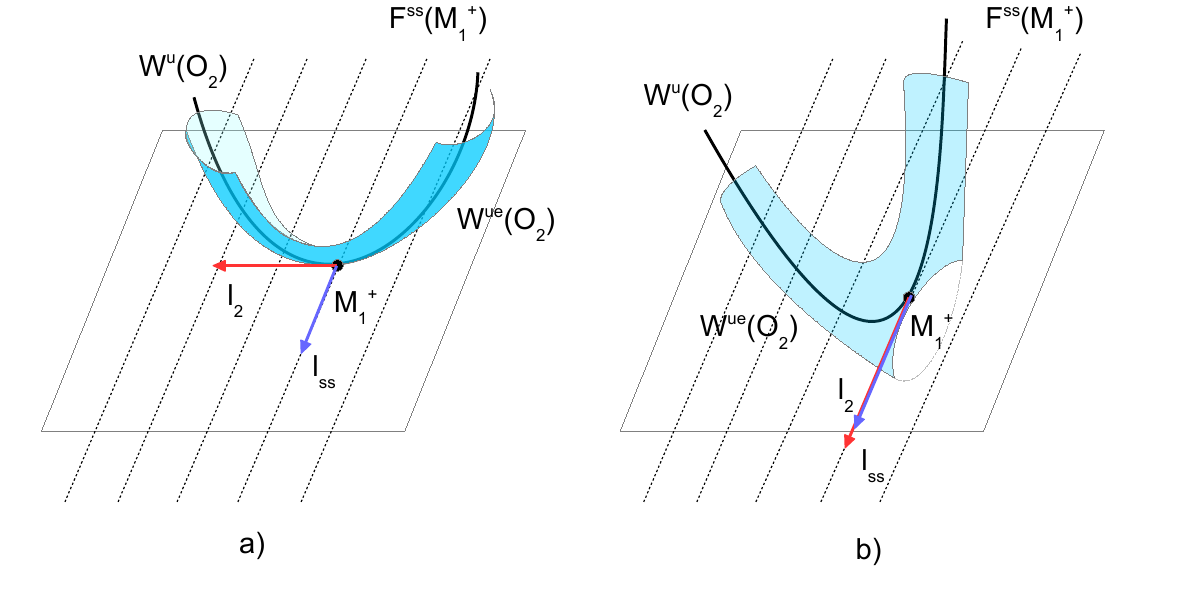}
 \caption{\label{fig05} Two types of the non-simple quadratic (heteroclinic) tangency:
(a) $W^{ue}(O_2)$ is transversal to $W^s_{loc}(O_1)$ and touches the leaf $F^{ss}(M_1^+)$; (b) $W^{ue}(O_2)$ is
tangent to $W^s_{loc}(O_1)$ and the curves $W^u(O_2)$ and $F^{ss}(M_1^+)$ have a general intersection at
$M_1^+$.} 
\end{figure}

Now consider diffeomorphism $f_0$. It has fixed points $O_1$ and $O_2$ and heteroclinic orbit $\Gamma_{21}$ in the points of which
manifolds $W^u(O_2)$ and $W^s(O_1)$ have a quadratic tangency. Everything mentioned above on the existence of extended unstable manifolds
and strong stable foliations can be applied to each of the saddles.
Let $U_{01} \ni O_1$ and $U_{02} \ni O_2$ be 
some small neighbourhoods of the fixed points, $M_1^+ \in W_{loc}^s(O_1) \subset U_{01}$ and 
$M_2^- \in W_{loc}^u(O_2) \subset U_{02}$ be two points of $\Gamma_{21}$
and $\Pi_1^+ \subset U_{01}$ and $\Pi_2^- \subset U_{02}$ their respective neighborhoods. 
Note that there exists some integer $q_1$ such that $M_1^+ = f_0^{q_1}(M_2^-)$.
Define the global map along $\Gamma_{21}$ for all small $\mu$ as $T_{21, \mu}: \Pi_2^- \to \Pi_1^+ = \left. f_\mu^{q_1} \right|_{\Pi_2^-}$ 
(for simpler notation, further we will omit the subscript $\mu$ for global and local maps, implicitly always assuming them to be 
the corresponding parametric families).
The heteroclinic tangency of $W^u(O_2)$ and $W^s(O_1)$ is called {\em simple} if image $T_{21}(P^{ue}(M_2^-))$
of tangent plane $P^{ue}(M_2^-)$ to $W^{ue}(O_2)$ at point $M_2^-$, 
intersects transversely the leaf $F_1^{ss}(M_1^+)$ of invariant foliation $F_1^{ss}$, containing point $M_1^+$.
If this condition is not fulfilled we call such a quadratic tangency {\em non-simple}. Following \cite{Tat01}, there
may be only two generic cases of non-simple heteroclinic tangencies defined by condition {\bf D1}:

Case ${\rm I}.\;$ {\it The surface $T_{21}(P^{ue}(M_2^-))$ is
transversal to the plane $W^{s}_{loc}(O_1)$  but is tangent to the line
$F_1^{ss}(M_1^+)$ at point $M_1^+$, Fig.\ref{fig05}~$($a$)$.}

Case ${\rm II}.\;$ {\it The surfaces $T_{21}(P^{ue}(M_2^-))$ and
$W^{s}_{loc}(O_1)$ have a quadratic tangency at $M_1^+$ and the curves $T_{21}(W^u_{loc}(O_2) \cap \Pi_2^-)$
and $F_1^{ss}(M_1^+)$ have a general intersection, Fig.\ref{fig05}~$($b$)$.}  \\

Thus, in Case ${\rm I}$ tangent vectors $l_u$ to
$T_{21}(W^u_{loc}(O_2))$ and $l_{ss}$ to $F_1^{ss}(M_1^+)$ are collinear,
while in Case ${\rm II}$ these vectors have different directions, the latter guarantees the absence of
additional degeneracies.

\begin{figure}[!ht]
\includegraphics[height=6cm]{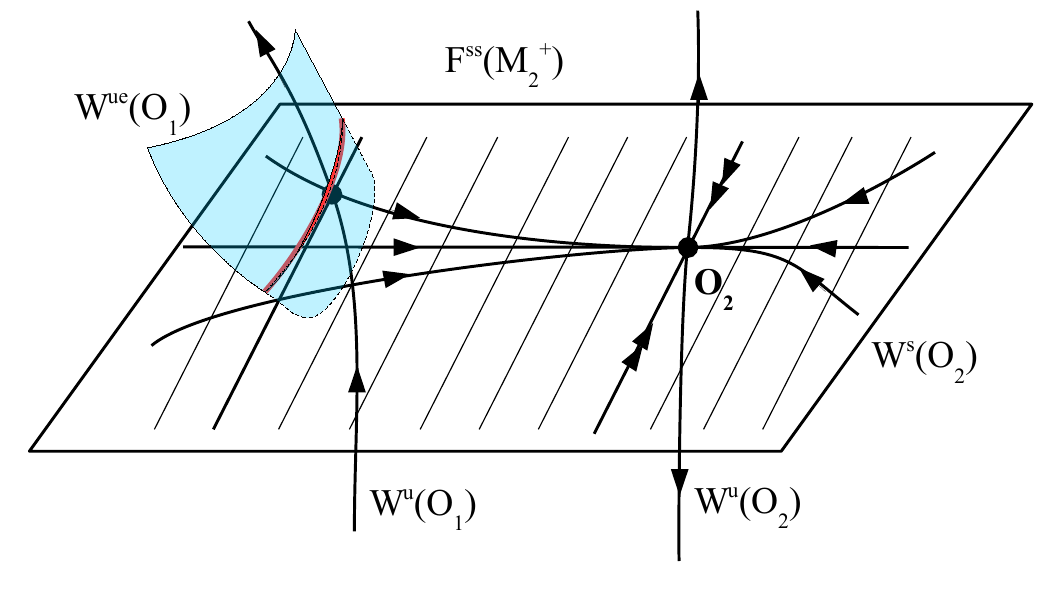} \caption{\label{Case3} A non-simple heteroclinic intersection of $W^u(O_1)$ and $W^s(O_2)$.}
\end{figure}

Note that if saddle points $O_1$ and $O_2$ coincide, we formally obtain the known definition of
a non-simple {\em homoclinic} tangency  \cite{Tat01, GOT14, GGT07}.
However, in distinct from the homoclinic case, the heteroclinic cycle under consideration allows
one more degeneracy, related to the second heteroclinic orbit $\Gamma_{12}$. This is the case when 
the transversal heteroclinic intersection of manifolds $W^u(O_1)$ and $W^s(O_2)$ is non-simple. 
Consider two heteroclinic points $M_1^- \in U_{01}$ and $M_2^+ \in U_{02}$ and their small respective neighbourhoods
$\Pi_1^- \subset U_{01}$ and $\Pi_2^+ \subset U_{02}$.
Again, there exists some integer $q_2$ such that $M_2^+ = f_0^{q_2}(M_1^-)$ so that
we define the global map from $U_{01}$ to $U_{02}$ as $T_{12}: \Pi_1^- \to \Pi_2^+ = \left. f_\mu^{q_2} \right|_{\Pi_1^-}$.
Let $P^{ue}(M_1^-)$ be the tangent plane to $W^{ue}(O_1)$ at $M_1^-$ and $F_2^{ss}(M_2^+)$ be the leaf of invariant foliation $F_2^{ss}$ on $W^s(O_2)$
passing through $M_2^+$. The heteroclinic intersection of
$W^u(O_1)$ and $W^s(O_2)$ is called {\em simple} if image $T_{12}(P^{ue}(M_1^-))$ and leaf $F_2^{ss}(M_2^+)$ 
intersect transversely. If this condition is not fulfilled the heteroclinic intersection is {\em non-simple}, see Fig.~\ref{Case3}. 
Thus, under condition {\bf D2}, we have

Case~${\rm III}.\;$ {\it The surface $T_{12}(P^{ue}(M_1^-))$ is
transversal to the plane $W^{s}_{loc}(O_2)$ but is tangent to the line
$F_2^{ss}(M_2^+)$ at $M_2^+$.}

In the present paper the birth of discrete Lorenz attractors in cases~${\rm I-III}$ is studied. The main result is given by the following 

\begin{teo}  \label{thmmain}
Let $f_\mu$ be the three-parametric family under consideration $(f_0$ satisfies  {\bf A}--{\bf D} and $f_\mu$ is a 
general unfolding of conditions {\bf B}, {\bf C} and {\bf D}, where {\bf D} is either {\bf D1} or {\bf D2}$)$. Then,
in any neighbourhood of the origin $\mu = 0$ in the parameter space there exist infinitely many domains
$\delta_{ij}$, where $\delta_{ij} \to (0,0,0)$ as $i,j \to \infty$, such that
the diffeomorphism $f_\mu$  has at $\mu \in \delta_{ij}$ a
discrete Lorenz attractor. \\
\end{teo}



\section{The first return map and the rescaling lemma}\label{res:het}

Let $U$ be a sufficiently small and fixed neighborhood of 
heteroclinic cycle $\{O_1, O_2, \Gamma_{12}, \Gamma_{21}\}$. It it composed 
as a union of small neighborhoods $U_{01}$ 
and $U_{02}$ of points $O_1$ and $O_2$ respectively, with a finite number of small neighborhoods $U_i$ of those
points of heteroclinic orbits $\Gamma_{12}$ and $\Gamma_{21}$
which do not belong to $U_{01} \cup U_{02}$. Each single-round periodic orbit lying entirely in
a small neighborhood of the heteroclinic cycle should have exactly one intersection point with each
of $U_i$ and the rest points lying in $U_{01} \cup U_{02}$.

Consider heteroclinic points $M_{1,2}^\pm$  and their respective small
neighborhoods $\Pi_{1,2}^\pm$ described in the previous section.
Define the first return map as a composition
of two local and two global maps. Local maps $T_{01}$ and $T_{02}$ are the restrictions of $f_\mu$ onto
$U_{01}$ and $U_{02}$ respectively and the global maps were defined in the following way:
$T_{12} = f_\mu^{q_1}: \Pi_1^- \to \Pi_2^+$, $T_{21} = f_\mu^{q_2}: \Pi_2^- \to \Pi_1^+$. 

Begin iterating $\Pi_1^+$ under the action of $T_{01}$. Starting from some number $i_0$ images
$T_{01}^{i} \Pi_1^+$ will have a nonempty intersection with $\Pi_1^-$. The same
applies for iterations of $\Pi_2^+$, there exists some $j_0$ such that for any $j \ge j_0$ 
image $T_{02}^{j} \Pi_2^+$ has a nonempty intersection with $\Pi_2^-$. 
The first return map $T_{ij} \equiv T_{21} T_{02}^j T_{12} T_{01}^i$ is thus defined on an infinite set of regions
$V_{ij}$ that lie in $\Pi_1^+$ and shrink to $M_1^+$ as $i, j \to \infty$. Their
images $f^i_\mu V_{ij}$ lie in $\Pi_1^-$, regions $f^{i + q_1}_\mu V_{ij}$ lie in $\Pi_2^+$, and  regions
$f^{i + q_1 + j}_\mu V_{ij}$ lie in $\Pi_2^-$, so $f^{i + q_1 + j + q_2}_\mu V_{ij}$ lie in
$\Pi_1^+$ again. 

To construct the first return map one needs first to write both local and global maps
in the most suitable form. For both local maps $T_{01}$ and $T_{02}$ it is the main normal form (\ref{t0norm}).
One its important property is that the iterations
$T_{0m}^k :\; U_{0m} \to U_{0m}$, $m = 1,2$, for any $k$ can be calculated in a
simple way, namely, in a form close to linear (see, for example, \cite{book, GS92}). 
Namely, for small $\mu$ iterations $T_{01}^k(\mu) : (x_0, y_0) \to (x_k, y_k)$ can be represented as:
\begin{equation}\label{eq:T0kk1}
\begin{array}{l}
  x_{1k} = \lambda_1^k x_{10} +
            \hat\lambda^{k} \xi_{1k}(x_0, y_k, \mu),\\
   x_{2k} = \hat\lambda^{k} \xi_{2k}(x_0, y_k, \mu),\\
  y_0  = \gamma_1^{-k} y_k + \hat \gamma_1^{-k} \xi_{3k}(x_0, y_k, \mu) ,
\end{array}
\end{equation}
and iterations $T_{02}^k(\mu): (u_0, v_0) \to (u_k, v_k)$ as
\begin{equation}\label{eq:T0kk2}
\begin{array}{l}
  u_{1k} = \nu_1^{k} u_{10} +
               \hat\nu^{k} \xi_{4k}(u_0, v_k, \mu),\\
  u_{2k} = \hat \nu^{k} \xi_{5k}(u_0, v_k, \mu),\\
  v_0  = \gamma_2^{-k} v_k + \hat\gamma_2^{-k} \xi_{6k}(u_0, v_k, \mu).
\end{array}
\end{equation}
Here $0 < \hat \lambda < \lambda_1, \;\; 0 < \hat \nu < \nu_1$, $\hat\gamma_{1,2} > \gamma_{1,2}$, 
functions $\xi_{mk}$ and their derivatives up to the order $(r - 2)$
are uniformly bounded and their higher order derivatives tend to zero, .

In normal coordinates (\ref{t0norm}) local stable and unstable
manifolds of $O_1$ in $U_1$ are $W^s_{loc} = \{ y = 0 \}$ and
$W^u_{loc} = \{ x = 0 \}$, the local
stable and unstable manifolds of $O_2$ in $U_2$ are $W^s_{loc} =
\{ v = 0 \}$ and $W^u_{loc} = \{ u = 0 \}$. Assume that for
$\mu = 0$, we have $M_1^- = (0, 0, y^-) \in W^u_{loc}(O_1)$, $M_2^+ = (u_1^+, u_2^+, 0)\in
W^s_{loc}(O_2)$,  and $M_2^- = (0, 0, v^-) \in W^u_{loc}(O_2)$,
$M_1^+ = (x_1^+, x_2^+, 0)\in W^s_{loc}(O_1)$.
Then  global maps  for all small $\mu$ are written as Taylor expansions near points $M_1^-$ and $M_2^-$:
\begin{equation}\label{eq:T12}
  T_{12} :
  \begin{array}{ccl}
    u - u^+ & = & A^{(1)} x + b^{(1)} (y - y^-) + O(\|x\|^2 +
                 \|x\| \cdot |y - y^-| + (y - y^-)^2), \\

    v & = & (c^{(1)})^{\top} x + d^{(1)} (y - y^-) + O(\|x\|^2 + \|x\|\cdot|y - y^-| + |y - y^-|^2),
  \end{array}
\end{equation}
\begin{equation}\label{eq:T21}
  T_{21} :
  \begin{array}{l}
    \bar x - x^+ = A^{(2)} u + b^{(2)} (v - v^-) + O(\|u\|^2 + \|u\|\cdot|v - v^-| + (v - v^-)^2), \\

    \bar y = y^+ + (c^{(2)})^{\top} u + d^{(2)} (v - v^-)^2 +  O(\|u\|^2 + \|u\|\cdot|v - v^-| + (v - v^-)^3),
  \end{array}
\end{equation}
where $d^{(1)} \neq 0$ and $d^{(2)} \neq 0$,
since $W^u(O_1)$ and $W^s(O_2)$ intersect transversely and the tangency between $W^u(O_2)$ and $W^s(O_1)$ is quadratic for $\mu = 0$.
Moreover, both maps $T_{12}$ and $T_{21}$ are diffeomorphisms, so that we have

\begin{equation}\label{eq:J12=}
J_{12} = \mbox{det}\; \begin{pmatrix}
a_{11}^{(1)} & a_{12}^{(1)} & b_1^{(1)} \\
a_{21}^{(1)} & a_{22}^{(1)} & b_2^{(1)} \\
c_1^{(1)} & c_2^{(1)} & d^{(1)}

\end{pmatrix}
\; \neq 0,\;\;\;
J_{21} = \mbox{det}\; \begin{pmatrix}
a_{11}^{(2)} & a_{12}^{(2)} & b_1^{(2)} \\
a_{21}^{(2)} & a_{22}^{(2)} & b_2^{(2)} \\
c_1^{(2)} & c_2^{(2)} & 0
\end{pmatrix}
\; \neq 0.
\end{equation}
In particular, this means that $\sqrt{b_1^{(2) 2} + b_2^{(2) 2}} \neq 0$ and
$\sqrt{c_1^{(2) 2} + c_2^{(2) 2}} \neq 0$ for $\mu = 0$. 

Now consider conditions \textbf{D1} and \textbf{D2} separately.

Case ${\rm I}$. Tangent plane $P^{ue}(M_2^-)$ to $W^{ue}_{loc}(O_2)$ at
 $M_2^-$ has equation ${u_2 = 0}$. The equation of
$T_{21}(P_{ue}(M_2^-))$ at $\mu = 0$ is obtained by putting $u_2 = 0$ into (\ref{eq:T21}).
Then the transversality of $T_{21}(P^{ue}(M_2^-))$ and $W^s_{loc}(O_1)$ which has the equation
$\bar y = 0$, yields $c_1^{(2)}(0) \neq 0$. The tangent vector to the
 line $T_{21}(P_{ue}(M_2^-)) \cap W^s_{loc}(O_1)$ at  point $M_1^+$ 
is $(b_1^{(2)}(0), b_2^{(2)}(0), 0)$. The equation of leaf $F^{ss}(M_1^+)$ is $\{x_1 = x_1^+, y = 0\}$.
Therefore, the tangency of $T_{21}(P^{ue}(M_2^-))$ and $F^{ss}(M_1^+)$
implies $b_1^{(2)}(0) = 0$. In this case $b_2^{(2)}(0) \neq 0$ and
$a_{11}^{(2)2} + a_{12}^{(2)2} \neq 0$ because of (\ref{eq:J12=}). 

Case $\rm II$. The equation of $T_1(P^{ue}(M_2^-))$ at
$\mu = 0$ is the same as in Case~$\rm I$. Then the tangency of surfaces
$T_{21}(P^{ue}(M_2^-))$ and $W^{s}_{loc}(O_1)$ at $\mu = 0$ implies
 $c_1^{(2)}(0) = 0$. Also, the tangent vectors to the lines
$T_{21}(W^u_{loc}(O_2) \cap \Pi_2^-)$ and $F^{ss}(M_1^+)$ at point $M_1^+$ are non-parallel if
$b_1^{(2)}(0) \neq 0$. 

Case $\rm III$. The equation of tangent plane $P^{ue}(M_1^-)$ to $W^{ue}_{loc}(O_1)$ at
$M_1^-$ is ${x_2 = 0}$. Putting this to (\ref{eq:T12}) gives the equation of its image
under $T_{12}$. The equation of leaf $F^{ss}(M_2^+)$ is $\{u_1 = u_1^+, v = 0\}$. Thus this leaf will
be tangent to $T_{12}(P^{ue}(M_1^-))$ at $\mu = 0$ if:
$$
A_{11}^{(1)}(0) = a_{11}^{(1)}(0) - b_1^{(1)}(0) c_1^{(1)}(0) / d_1(0) = 0.
$$


We are now able to write the non-simple heteroclinic orbit conditions in the explicit form for all three cases:

\begin{equation}
\begin{array}{c}
 {\rm Case \; I}: \: b_1^{(2)}(0) = 0, \; b_2^{(2)}(0) \neq 0, \; c_1^{(2)}(0) \neq 0, \; A_{11}^{(1)}(0) \neq 0. \\

 {\rm Case \; II}: \: b_1^{(2)}(0) \neq 0, \; c_1^{(2)}(0) = 0, \; c_2^{(2)}(0) \neq 0, \; A_{11}^{(1)}(0) \neq 0. \\
 
 {\rm Case \; III}: \: b_1^{(2)}(0) \neq 0, \; c_1^{(2)}(0) \neq 0, \; A_{11}^{(1)}(0) = 0.
 \label{eq:deg_cond}
\end{array}
\end{equation}

We will construct a three-parameters family $f_\mu$, where $\mu = (\mu_1, \mu_2, \mu_3)$,
As the first governing parameter we take
the splitting parameter $\mu_1$ for the quadratic heteroclinic tangency so that
\begin{equation}
\;\;\mu_1 \equiv y^+(\mu)\;.
\label{mu1het}
\end{equation}
The second parameter $\mu_2$ is responsible for the degeneracy related to, respectively, conditions ({\bf D1}):
\begin{equation}
\mu_2 = b_1^{(2)}\;\;\;{\rm in\; Case\; I},
\label{mu2Ihet}
\end{equation}
\begin{equation}
\mu_2 = c_1^{(2)}\;\;\;
{\rm in\; Case\; II}
\label{mu2IIhet}
\end{equation}
or ({\bf D2}):
\begin{equation}
\displaystyle \mu_2 = A_{11}^{(1)} = a_{11}^{(1)} - \frac{b_1^{(1)} c_{1}^{(1)}}{d_1}\;\;\;
{\rm in\; Case\; III.}
\label{mu2IIIhet}
\end{equation}


The third parameter has been already given by formula (\ref{eq:e3}).

\begin{lem}
Let $f_{\mu_1, \mu_2, \mu_3}$ be the family under consideration. Then, in the
space $(\mu_1, \mu_2, \mu_3)$ there exist infinitely many regions $\Delta_{ij}$
ac\-cu\-mu\-la\-ting to the origin as $i, j \to \infty$, such that the first return map $T_{ij}$ in appropriate rescaled
coordinates and parameters is asymptotically $C^{r - 1}$-close to one of the following limit maps.

{\rm 1)} In Case I, the limit map is
\begin{equation}
\label{HIhet}
\begin{array}{l}
\bar X_1 \; = \; -J_{ij} X_2 + M_2 Y,\;\; \bar X_2 \; = \; Y,\;\;
\bar Y = M_1 - X_1 - Y^2,
\end{array}
\end{equation}
where
\begin{equation}
\label{M123Ihet}
\begin{array}{l}
M_1 = -d^{(1)^2} d^{(2)} \gamma_1^{2i} \gamma_2^{2j}
     (\mu_1 + \nu_1^j c_1^{(2)} u_1^+ - \gamma_1^{-i}y^- + \nu_{ij}^1),\\
M_2 = (\mu_2 + \rho_{ij}^1)c_1^{(2)} A_{11}^{(1)} \lambda_1^i \gamma_1^i \nu_1^j \gamma_2^j, \\
J_{ij} = J_{12}J_{21}(\lambda_1\lambda_2\gamma_1)^i(\nu_1\nu_2\gamma_2)^j,
\end{array}
\end{equation}
and $\nu_{ij}^1 = O(\hat\gamma_1^{-i} + \hat\nu^j + \gamma_1^{-i}\gamma_2^{-j})$, $\rho_{ij}^1 = O(\nu_1^j)$.

{\rm 2)}  In Case II, the limit map is
\begin{equation}
\label{HIIhet}
\begin{array}{l}
\bar X_1 \; = \; Y,\;\; \bar X_2 \; = \; X_1,\;\;
\bar Y = M_1 + M_2 X_1 + B X_2- Y^2,
\end{array}
\end{equation}
where
\begin{equation}
\label{M123IIhet}
\begin{array}{l}
M_1 = -d^{(1)^2} d^{(2)} \gamma_1^{2i} \gamma_2^{2j}
     (\mu_1 + \nu_1^j \mu_2 u_1^+ - \gamma_1^{-i}y^- + \nu_{ij}^2),\\

M_2 = (\mu_2 + \rho_{ij}^2 )b_1^{(2)} A_{11}^{(1)} \lambda_1^i \nu_1^j \gamma_1^{i} \gamma_2^j, \\
B = J_{12}J_{21}(\lambda_1\lambda_2\gamma_1)^i(\nu_1\nu_2\gamma_2)^j
\end{array}
\end{equation}
and $\nu_{ij}^2 = O(\hat\gamma_1^{-i} + \hat\nu^j + \gamma_1^{-i}\gamma_2^{-j} + \nu_1^j \left( \hat\lambda / \lambda_1 \right)^i )$,
$\rho_{ij}^2 = O\left( \left( \hat\lambda / \lambda_1 \right)^i + \left( \hat\nu / \nu_1 \right)^j \right)$.

{\rm 3)}  In Case III, the limit map is
\begin{equation}
\label{HIIIhet}
\begin{array}{l}
\bar X_1 \; = \; Y,\;\; \bar X_2 \; = \; X_1,\;\;
\bar Y = M_1 + M_2 X_1 + B X_2- Y^2,
\end{array}
\end{equation}
where
\begin{equation}
\label{M123IIIhet}
\begin{array}{l}
M_1 = -d^{(1)^2} d^{(2)} \gamma_1^{2i} \gamma_2^{2j}
      (\mu_1 + \nu_1^j c_1^{(2)} u_1^+ - \gamma_1^{-i}y^- + \nu_{ij}^3),\\

M_2 = (\mu_2 + \rho_{ij}^3 )b_1^{(2)} c_1^{(2)} \lambda_1^i \nu_1^j \gamma_1^{i} \gamma_2^j, \\
B = J_{12}J_{21}(\lambda_1\lambda_2\gamma_1)^i(\nu_1\nu_2\gamma_2)^j
\end{array}
\end{equation}
and $\nu_{ij}^3 = O(\hat\gamma_1^{-i} + \hat\nu^j + \gamma_1^{-i}\gamma_2^{-j})$, 
$\rho_{ij}^3 = O\left( \left(\hat\lambda / \lambda_1\right)^{i} + \left(\hat\nu/\nu_1\right)^j \right)$.
\label{th4}
\end{lem}

It is easy to see that the rescaled first return map in Cases~II and~III is asymptotically close to the 3D Henon map (\ref{H3D}). In
Case~I for system (\ref{HIhet}) we make an additional change of coordinates $X_{1new} = X_1 - M_2 X_2$ and scale $X_1$ by $(-J_{ij})$, bringing it
again to the form (\ref{H3D}). Next, the statement of Theorem~\ref{thmmain} follows immediately -- as shown in \cite{GOST05, GMO06}, this 
three-dimensional Henon map
possesses the discrete Lorenz attractor for an open set of parameters $(M_1, M_2, B)$. Hence for each sufficiently large $i$ and $j$, for which
the Jacobian $J_{ij}$ stays finite, the corresponding domain $\delta_{ij}$ in the original parameters $(\mu_1, \mu_2, \mu_3)$ is determined from
formulas (\ref{M123Ihet}), (\ref{M123IIhet}) or (\ref{M123IIIhet}) respectively in Cases~I--III. These domains accumulate to the origin when
$i$ and $j$ unboundedly grow. This proves Theorem~\ref{thmmain}.



\section{Proof of the rescaling Lemma~\ref{th4}.}\label{sec:th3proof}

Note that the first return map $T_{ij}$ is rescaled differently 
in cases~$\rm I$--$\rm III$, however, there is a preparation part of the proof that is conducted in the same way 
for all the cases.


Using formulas (\ref{eq:T0kk1}), (\ref{eq:T0kk2}), (\ref{eq:T12}) and (\ref{eq:T21}), we obtain the following expression for the first
return map $ T_{ij}\equiv T_{21} T_{02}^jT_{12} T_{01}^i: \Pi^+_1 \to\Pi^+_1$
\begin{equation}\label{eq:Tij12}
    \begin{array}{l}
    u_1 - u_1^+ = a_{11}^{(1)} (\lambda_1^i x_1 + \hat\lambda^{i} \xi_{1i}(x, y, \mu)) +    
            a_{12}^{(1)} \hat \lambda^{i} \xi_{2i}(x, y, \mu) +  b_1^{(1)} (y - y^-) + \\  
            \qquad + O(\lambda_1^{2i}\|x\|^2 + \lambda_1^i \|x\| \cdot |y - y^-| + (y - y^-)^2), \\
    
    u_2 - u_2^+ = a_{21}^{(1)} (\lambda_1^i x_1 + \hat\lambda^{i}\xi_{1i}(x, y, \mu)) +    
            a_{22}^{(1)} \hat\lambda^{i} \xi_{2i}(x, y, \mu) + b_2^{(1)} (y - y^-) +  \\ 
            \qquad  + O(\lambda_1^{2i}\|x\|^2 + \lambda_1^i \|x\| \cdot |y - y^-| + (y - y^-)^2), \\

    \gamma_2^{-j}(v + \left( \hat\gamma_2 / \gamma_2 \right)^{-j} \xi_{6j}(u, v, \mu))) = 
    c_1^{(1)} (\lambda_1^i x_1 + \hat\lambda^{i} \xi_{1i}(x, y, \mu)) + 
    c_2^{(1)} \hat\lambda_1^{i} \xi_{2i}(x, y, \mu) + \\
      \qquad + d^{(1)}(y - y^-) + O(\lambda_1^{2i}\|x\|^2 + \lambda_1^i \|x\| \cdot |y - y^-| + (y - y^-)^2),
     \end{array}
\end{equation}
\begin{equation}\label{eq:Tij121}     
     \begin{array}{l}
    \bar x_1 - x_1^+ = a_{11}^{(2)} (\nu_1^j u_1 + \hat\nu^{j} \xi_{4j}(u, v, \mu)) +
       a_{12}^{(2)} \hat\nu^j  \xi_{5j}(u, v, \mu) 
        + b_1^{(2)} (v - v^-) + \\ 
        \qquad  + O(\nu_1^{2j} \|u\|^2 + \nu_1^j |u| \cdot |v - v^-| + (v - v^-)^2), \\

    \bar x_2 - x_2^+ = a_{21}^{(2)} (\nu_1^j u_1 + \hat\nu^{j} \xi_{4j}(u, v, \mu)) +
       a_{22}^{(2)} \hat\nu^j \xi_{5j}(u, v, \mu) + b_2^{(2)} (v - v^-) + \\
      \qquad + O(\nu_1^{2j} \|u\|^2 + \nu_1^j |u| \cdot |v - v^-| + (v - v^-)^2), \\

    \gamma_1^{-i} (\bar y + \left( \hat\gamma_1 / \gamma_1 \right)^{-i} \xi_{3i}(\bar x, \bar y, \mu))) = \mu_1 + 
    c_1^{(2)} (\nu_1^j u_1 + \hat\nu^j \xi_{4j}(u, v, \mu)) + \\ 
    \qquad + c_2^{(2)} \hat\nu^j \xi_{5j}(u, v, \mu))
        + d^{(2)} (v - v^-)^2 + O(\nu_1^{2j} \|u\|^2 + \nu_1^j \|u\| \cdot |v - v^-| + \\ 
        \qquad + (v - v^-)^3),
  \end{array}
\end{equation}

Make a coordinate shift
$u_{new} = u - u^+ + \varphi_{ij}^1$, $v_{new} = v - v^- +
\varphi_{ij}^2$, $x_{new} = x - x^+ + \psi_{ij}^1$, $y_{new} = y - y^- + \psi_{ij}^2$,
where $\varphi_{ij}^1,\; \psi_{ij}^2 = O(\gamma_2^{-j} + \lambda_1^i)$ and 
$\varphi_{ij}^2,\;  \psi_{ij}^1 = O(\nu_1^j)$. With that, the nonlinearity functions in the left parts of the third equations 
of (\ref{eq:Tij12}) and (\ref{eq:Tij121}) can be expressed as Taylor expansions 
$\xi_{6j}(u + u^+ + \varphi_{ij}^1, v + v^- + \varphi_{ij}^2, \mu)) = 
\xi_{6j}^0 + \xi_{6j}^1 v + \xi_{6j}^2(u, v) + \xi_{6j}^3(v)$, $\xi_{3i}(\bar x + x^+ + \psi_{ij}^1, \bar y + y^- + \psi_{ij}^2, \mu)) = 
\xi_{3i}^0 + \xi_{3i}^1 \bar y + \xi_{3i}^2(\bar x, \bar y) + \xi_{3i}^3(\bar y)$ respectively, where coefficients  
$\xi_{6j}^0$, $\xi_{6j}^1$, $\xi_{3i}^0$, $\xi_{3i}^1$ are uniformly bounded in $i$ and $j$ for all small $\mu$ and
$\xi_{6j}^2 (u, v) = O(u)$, $\xi_{3i}^2 (\bar x, \bar y) = O(\bar x)$, $\xi_{6j}^3(v) = O(v^2)$, $\xi_{3i}^3(\bar y) = O(\bar y^2)$. 
We select constants $\varphi_{ij}^1, \; \varphi_{ij}^2,\; \psi_{ij}^1,\; \psi_{ij}^2$ in such a way
that all constant terms in equations (\ref{eq:Tij12}), the constant terms in the first two equations and the linear in $v_{new}$ term
in the last equation of (\ref{eq:Tij121}) vanish. In addition, we plug the expressions for $u$ coordinates from the
first two equations of (\ref{eq:Tij12}) into the third one, this will cause additions of order $\hat\gamma_2^{-j}$ to all
the coefficients. The system is rewritten as:
\begin{equation}\label{eq:3}
    \begin{array}{l}
    u_1 = a_{11}^{(1)} \lambda_1^i x_1 + O(\hat\lambda^i + \lambda_1^i \gamma_2^{-j}) O(\|x\|) + b_1^{(1)} y +   
            \lambda_1^i O(\|x\| \cdot |y|) + O(y^2), \\
    
    u_2 = a_{21}^{(1)} \lambda_1^i x_1 + O(\hat\lambda^i + \lambda_1^i \gamma_2^{-j}) O(\|x\|) + b_2^{(1)} y + 
            \lambda_1^i O(\|x\| \cdot |y|) + O(y^2), \\

    \gamma_2^{-j}(1 + q_{ij}^{(2)})v + \hat\gamma_2^{-j} O(v^2) = 
    c_1^{(1)} \lambda_1^i x_1 + O(\hat\lambda^i + \lambda_1^i \gamma_2^{-j}) O(\|x\|) + d^{(1)} y + \\
         \qquad + \lambda_1^i O(\|x\| \cdot |y|) + O(y^2),
     \end{array}
\end{equation}
\begin{equation}\label{eq:4}     
     \begin{array}{l}
    \bar x_1 = a_{11}^{(2)} \nu_1^j u_1 + \tilde a_{12}^{(2)} \hat\nu^j u_2 + b_1^{(2)} v + 
         O(\hat\nu^j \|u\|^2 + \nu_1^j \|u\| \cdot |v| + v^2), \\

    \bar x_2 = a_{21}^{(2)} \nu_1^j u_1 + \tilde a_{22}^{(2)} \hat\nu^j u_2 + b_2^{(2)} v + 
         O(\hat\nu^j \|u\|^2 + \nu_1^j \|u\| \cdot |v| + v^2), \\

    \gamma_1^{-i} (1 + q_{ij}^{(1)})\bar y + \hat\gamma_1^{-i} O(\bar x) + \hat\gamma_1^{-i} O(\bar y^2) = M^1 + c_1^{(2)} \nu_1^j u_1 + 
    \tilde c_2^{(2)} \hat\nu^j u_2 + d^{(2)} v^2 + \\ \qquad + O(\hat\nu^{j} \|u\|^2 + \nu_1^j \|u\| \cdot |v| + |v|^3),
  \end{array}
\end{equation}
where $q_{ij}^{(1)} = O\left(\left( \hat\gamma_1 / \gamma_1 \right)^{-i} \right)$, 
$q_{ij}^{(2)} = O\left(\left( \hat\gamma_2 / \gamma_2 \right)^{-j} \right)$, coefficients marked with ``tilde'' are uniformly bounded
for small $\mu$ and the following expression is valid for $M^1$:

\begin{equation} \label{MfFeps}
M^1 = \mu_1 + \nu_1^j c_1^{(2)} u_1^+ - \gamma_1^{-i}y^- + O(\hat\gamma_1^{-i} + \hat\nu^j + \gamma_1^{-i}\gamma_2^{-j}). 
\end{equation}

Next, we take the right-hand side of the third equation of (\ref{eq:3}) divided by factor 
$\gamma_2^{-j}(1 + q_{ij}^{(2)})$ from the left-hand side as the new variable $y$ --  the equation becomes the following
$v + \left(\left( \hat\gamma_2 / \gamma_2 \right)^{-j} \right) O(v^2) = y$. We substitute this formula instead of the $y$ variable 
to all equations. Defining $u_{new} = u - (b^{(1)} / d^{(1)}) \gamma_2^{-j} v + O(\hat\gamma_2^{-j} v^2)$ we
eliminate all terms in the equation for $u$ which depend
on $v$ alone. In addition, we substitute the expressions for $\bar x$ to the last equation of (\ref{eq:4}). These actions cause
the linear in $v$ term of order $O(\hat\gamma_1^{-i} + \nu_1^j \gamma_2^{-j})$ to appear in the equation for $\bar v$. We will
make it zero again later.
Thus, we obtain
\begin{equation}\label{eq:Tijshc}
    \begin{array}{l}
 u_1 = A_{11}^{(1)} \lambda_1^i x_1 + O(\hat\lambda^i + \lambda_1^i \gamma_2^{-j}) O(\|x\|) +  \lambda_1^i\gamma_2^{-j}O( \|x\| \cdot |v|), \\
 
 u_2 = A_{21}^{(1)} \lambda_1^i x_1 +O(\hat\lambda^i + \lambda_1^i \gamma_2^{-j}) O(\|x\|) + \lambda_1^i\gamma_2^{-j}O( \|x\| \cdot |v|), \\

\bar x_1 = a_{11}^{(2)} \nu_1^j u_1 + \tilde a_{12}^{(2)} \hat\nu^j u_2 + b_1^{(2)} v +
    O(\hat\nu^j \|u\|^2 + \nu_1^j\|u\| \cdot |v| + v^2), \\

\bar x_2 = a_{21}^{(2)} \nu_1^j u_1 + \tilde a_{22}^{(2)} \hat\nu^j u_2 + b_2^{(2)} v +
    O(\hat\nu^j \|u\|^2 + \nu_1^{j}\|u\| \cdot |v| + v^2), \\

\displaystyle  \frac{\gamma_1^{-i}\gamma_2^{-j}}{d^{(1)}} \bar v(1 + q_{ij}^{(3)}) + \gamma_1^{-i}\hat\gamma_2^{-j} O(\bar v^2) =  M^1 + c_1^{(2)} \nu_1^j u_1 + 
    \tilde c_2^{(2)} \hat\nu^j u_2 + \\ 
      \qquad\qquad\qquad\qquad  + O(\hat\gamma_1^{-i} + \nu_1^j \gamma_2^{-j}) v 
      + d^{(2)} v^2 + O(\hat\nu^j \|u\|^2 + \nu_1^{j}\|u\| \cdot |v| + |v|^3),
  \end{array}
\end{equation}
where $q_{ij}^{(3)} = O\left(\left( \hat\gamma_1 / \gamma_1 \right)^{-i} + \left(\hat\gamma_2 / \gamma_2 \right)^{-j} \right)$ and 
\begin{equation} \label{eq:Anew}
A_{11}^{(1)} = a_{11} - b_1^{(1)} c_1^{(1)} / d^{(1)}, \; A_{21}^{(1)} = a_{21} - b_2^{(1)} c_1^{(1)} / d^{(1)}.
\end{equation}

Next, we  substitute $u$ as a function of $x$ and $v$  from the first two equations to the last three ones. After this, in addition, 
we make a shift of $(x, v)$ coordinates by a constant of order $O(\hat\gamma_1^{-i} + \nu_1^j \gamma_2^{-j})$ to 
nullify the linear in $v$ term in the last equation. This gives
us the following formula for the map $T_{ij}: (x, v) \mapsto (\bar x, \bar v)$:
\begin{equation}\label{eq:Tij4}
  \begin{array}{l}
    \bar x_1 = A_{11}^{(1)} a_{11}^{(2)} \lambda_1^i \nu_1^j x_1 + \tilde a_{12} s_{ij}^{(1)} x_2 + 
    b_1^{(2)} v + O(\hat\lambda^i \nu_1^j + \lambda_1^i \nu_1^j \gamma_2^{-j}) O(\|x\|^2) + \\ 
    \qquad\qquad + \lambda_1^i \nu_1^j O(\|x\| \cdot |v|) + O(v^2), \\

    \bar x_2 = A_{11}^{(1)} a_{21}^{(2)} \lambda_1^i \nu_1^j x_1 + \tilde a_{22} s_{ij}^{(2)} x_2 + 
    b_2^{(2)} v + O(\hat\lambda^i \nu_1^j + \lambda_1^i \nu_1^j \gamma_2^{-j}) O(\|x\|^2) + \\ 
    \qquad\qquad + \lambda_1^i \nu_1^j O(\|x\| \cdot |v|) + O(v^2), \\

    \displaystyle  \frac{\gamma_1^{-i}\gamma_2^{-j}}{d^{(1)}} \bar v(1 + q_{ij}^{(3)}) + \gamma_1^{-i}\hat\gamma_2^{-j} O(\bar v^2) = 
    M^1 + A_{11}^{(1)} c_1^{(2)} \lambda_1^i \nu_1^j  x_1 +
    \tilde c_{2} s_{ij}^{(3)} x_2 + d^{(2)} v^2 + \\

    \qquad\qquad + O(\hat\lambda^i \nu_1^j + \lambda_1^i \nu_1^j \gamma_2^{-j}) O(\|x\|^2) + \lambda_1^i \nu_1^j O(\|x\| \cdot |v|) + O(v^3),
 \end{array}
\end{equation}
where $s_{ij}^{(k)} = O(\hat\lambda^i \nu_1^j + \lambda_1^i \nu_1^j \gamma_2^{-j})$.


{\bf Case~{\rm I}}. The second parameter in the first case is introduced as $\mu_2\equiv b_1^{(2)}(\mu)$ and we also recall that $c_1^{(2)}$, $b_2^{(2)}$
and $A_{11}^{(1)}$ are bounded from zero due to (\ref{eq:deg_cond}). We make the following change of coordinates:
$$
\displaystyle
 x_{1new} = x_{1} + O\left(\left( \hat\lambda/\lambda_1\right)^{i} + \gamma_2^{-j} \right) O(\| x \|)\;, \;
x_{2new} = x_{2}\;,\; v_{new} = v
$$
such that all the terms which depend only on $x$--coordinates are now put into $x_{1new}$ in the third equation. 
Then (\ref{eq:Tij4}) is rewritten in the form:
\begin{equation}
\label{tk-1+-}
\begin{array}{l}
 \bar x_1 = A_{11}^{(1)} a_{11}^{(2)} \lambda_1^i \nu_1^j x_1 + \tilde a_{12} s_{ij}^{(1)} x_2 + 
    (\mu_2 + \rho_{ij}^1) v + O(\hat\lambda^i \nu_1^j + \lambda_1^i \nu_1^j \gamma_2^{-j}) O(\|x\|^2) + \\ 
    \qquad\qquad + \lambda_1^i \nu_1^j O(\|x\| \cdot |v|) + O(v^2), \\

    \bar x_2 = A_{11}^{(1)} a_{21}^{(2)} \lambda_1^i \nu_1^j x_1 + \tilde a_{22} s_{ij}^{(2)} x_2 + 
    b_2^{(2)} v + O(\hat\lambda^i \nu_1^j + \lambda_1^i \nu_1^j \gamma_2^{-j}) O(\|x\|^2) + \\ 
    \qquad\qquad + \lambda_1^i \nu_1^j O(\|x\| \cdot |v|) + O(v^2), \\

    \displaystyle  \frac{\gamma_1^{-i}\gamma_2^{-j}}{d^{(1)}} \bar v(1 + q_{ij}^{(3)}) + \gamma_1^{-i}\hat\gamma_2^{-j} O(\bar v^2) = 
    M^1 + A_{11}^{(1)} c_1^{(2)} \lambda_1^i \nu_1^j  x_1 + d^{(2)} v^2 + \\

    \qquad\qquad + \lambda_1^i \nu_1^j O(\|x\| \cdot |v|) + O(v^3),

\end{array}
\end{equation}
where $\rho_k^1 = O(\nu_1^{j})$.
Now we rescale the coordinates as follows
$$
\displaystyle v  = -\frac{\gamma_1^{-i}\gamma_2^{-j}}{d^{(1)} d^{(2)}} (1 + q_{ij}^{(3)}) \;Y \;,\; 
x_1  = \frac{\lambda_1^{-i} \gamma_1^{-2i} \nu_1^{-j} \gamma_2^{-2j}}
{c_1^{(2)} A_{11}^{(1)} (d^{(1)})^2 d^{(2)}} \;X_1\;,\;
x_2  = - \frac{b_2^{(2)}\gamma_1^{-i}\gamma_2^{-j}} {d^{(1)} d^{(2)}} \;X_2.
$$
Then system (\ref{tk-1+-}) is rewritten in the new coordinates:
\begin{equation}
\label{tkres1}
\begin{array}{l}
\displaystyle \bar X_1 \; = \; - J_{ij} X_2 +
M_2 Y  + O(\lambda_1^i \nu_1^j)\;,\; \\
\bar X_2  \; = \; Y + O(\hat\lambda^i \nu_1^j + \lambda_1^i \nu_1^j \gamma_2^{-j} + \gamma_1^{-i} \gamma_2^{-j})\;,\; \\
\bar Y \; = \; M_1 - X_1 - Y^2 + O(\gamma_1^{-i} \gamma_2^{-j})\;,
\end{array}
\end{equation}
where formulas (\ref{M123Ihet}) are valid for $M_1$ and $M_2$. Note, that the coefficient $J_{ij}$
is the Jacobian of map (\ref{tkres1}) up to asymptotically small in $i,j$ terms,
and, hence, $J_{ij}$ coincides in the main order with the Jacobian of map $T_{21} T_{02}^jT_{12} T_{01}^i$, i.e. it is given by formula (\ref{M123Ihet}).



{\bf Case~{\rm II}}. Now we have $\mu_2 \equiv c_1^{(2)}(\mu)$ and coefficients $b_1^{(2)}$, $c_2^{(2)}$ and $A_{11}^{(1)}$ are not zeros.  
Introduce the new coordinates as $ x_{1new} = x_{1} \;,\; x_{2new} = x_{2} - (b_2^{(2)} / b_1^{(2)})
x_1 \;,\; v_{new} = v$. Then (\ref{eq:Tij4}) recasts as
\begin{equation}
\label{tk-1+n}
\begin{array}{l}
\bar x_1 = A_{11}^{(1)} a_{11}^{(2)} \lambda_1^i \nu_1^j x_1 + \tilde a_{12} s_{ij}^{(1)} x_2 + 
    b_1^{(2)} v + O(\hat\lambda^i \nu_1^j + \lambda_1^i \nu_1^j \gamma_2^{-j}) O(\|x\|^2) + \\ 
    \qquad\qquad + \lambda_1^i \nu_1^j O(\|x\| \cdot |v|) + O(v^2), \\

    \bar x_2 = A_{11}^{(1)} A_{21}^{(2)} \lambda_1^i \nu_1^j x_1 + \tilde a_{22} s_{ij}^{(2)} x_2 + 
     O(\hat\lambda^i \nu_1^j + \lambda_1^i \nu_1^j \gamma_2^{-j}) O(\|x\|^2) + \\ 
    \qquad\qquad + \lambda_1^i \nu_1^j O(\|x\| \cdot |v|) + O(v^2), \\

    \displaystyle  \frac{\gamma_1^{-i}\gamma_2^{-j}}{d^{(1)}} \bar v(1 + q_{ij}^{(3)}) + \gamma_1^{-i}\hat\gamma_2^{-j} O(\bar v^2) = 
    M^1 + (\mu_2 + \rho_{ij}^2) A_{11}^{(1)} \lambda_1^i \nu_1^j  x_1 + \tilde c_{2} s_{ij}^{(3)} x_2 + \\ 
    
    \qquad\qquad + d^{(2)} v^2 + O(\hat\lambda^i \nu_1^j + \lambda_1^i \nu_1^j \gamma_2^{-j}) O(\|x\|^2) + \lambda_1^i \nu_1^j O(\|x\| \cdot |v|) + O(v^3),

\end{array}
\end{equation}
where $\rho_{ij}^2 = O\left(\left(\hat\lambda / \lambda_1 \right)^{i} + \left(\hat\nu / \nu_1\right)^{j}\right)$, 
$A_{21}^{(2)} = a_{21}^{(2)} - (b_2^{(2)} / b_1^{(2)}) a_{11}^{(2)} \neq 0$ due to (\ref{eq:J12=}) and (\ref{eq:deg_cond}).  
Now we rescale the coordinates as follows
$$
\displaystyle  v = -\frac{\gamma_1^{-i}\gamma_2^{-j}}{d^{(1)} d^{(2)}} (1 + q_{ij}^{(3)}) \; Y \;,\; 
x_1  = - \frac{b_1^{(2)} \gamma_1^{-i} \gamma_2^{-j}} {d^{(1)} d^{(2)}}\;X_1 \;,\; 
x_2  = - \frac{b_1^{(2)} A_{11}^{(1)} A_{21}^{(2)} \lambda_1^i \gamma_1^{-i} \lambda_2^j \gamma_2^{-j}}{d^{(1)} d^{(2)}}\;X_2 \;.
$$
After this, we can rewrite (\ref{tk-1+n}) in the following form
\begin{equation}
\label{tkres21}
\begin{array}{l}
 \bar X_1 \; = \;  Y + O( \lambda_1^i \nu_1^j  + \gamma_1^{-i} \gamma_2^{-j}))\;, \\
\bar X_2  \; = \;  X_{1} + O( \lambda_1^i \nu_1^j  + \gamma_1^{-i} \gamma_2^{-j})) \;,\\
 \bar Y = M_1 + M_2 X_1 + J_{ij} X_2 - Y^2 + O( \lambda_1^i \nu_1^j  + \gamma_1^{-i} \gamma_2^{-j})).
\\
\end{array}
\end{equation}
and formulas (\ref{M123IIhet}) are valid for $M_1$, $M_2$ and $J_{ij}$.


{\bf Case~{\rm III}}. Here we have $\mu_2 \equiv A_{11}^{(1)}  = a_{11} - b_1^{(1)} c_1^{(1)} / d^{(1)}$, 
and coefficients $b_1^{(2)}$, $c_1^{(2)}$ are not zeros.  
Introduce the new coordinates as in the previous case: $ x_{1new} = x_{1} \;,\; x_{2new} = x_{2} - (b_2^{(2)} / b_1^{(2)})
x_1 \;,\; v_{new} = v$. The system (\ref{eq:Tij4}) is then rewritten as:
\begin{equation}\label{eq:case3}
\begin{array}{l}
   \bar x_1 =  a_{11}^{(2)} (\mu_2 + \rho_{ij}^4) \lambda_1^i \nu_1^j x_1 + \tilde a_{12} s_{ij}^{(1)} x_2 + 
    b_1^{(2)} v + O(\hat\lambda^i \nu_1^j + \lambda_1^i \nu_1^j \gamma_2^{-j}) O(\|x\|^2) + \\ 
    \qquad\qquad + \lambda_1^i \nu_1^j O(\|x\| \cdot |v|) + O(v^2), \\

    \bar x_2 = A_{21}^{(2)} (\mu_2 + \rho_{ij}^5) \lambda_1^i \nu_1^j x_1 + \tilde a_{22} s_{ij}^{(2)} x_2 + 
     O(\hat\lambda^i \nu_1^j + \lambda_1^i \nu_1^j \gamma_2^{-j}) O(\|x\|^2) + \\ 
    \qquad\qquad + \lambda_1^i \nu_1^j O(\|x\| \cdot |v|) + O(v^2), \\

    \displaystyle  \frac{\gamma_1^{-i}\gamma_2^{-j}}{d^{(1)}} \bar v(1 + q_{ij}^{(3)}) + \gamma_1^{-i}\hat\gamma_2^{-j} O(\bar v^2) = 
    M^1 + c_1^{(2)} (\mu_2 + \rho_{ij}^3)\lambda_1^i \nu_1^j  x_1 + \tilde c_{2} s_{ij}^{(3)} x_2 + \\ 
    
    \qquad\qquad + d^{(2)} v^2 + O(\hat\lambda^i \nu_1^j + \lambda_1^i \nu_1^j \gamma_2^{-j}) O(\|x\|^2) + \lambda_1^i \nu_1^j O(\|x\| \cdot |v|) + O(v^3),
\end{array}
\end{equation}
where $A_{21}^{(2)} = a_{21}^{(2)} - (b_2^{(2)} / b_1^{(2)}) a_{11}^{(2)}$ and
$\rho_{ij}^{3,4,5} = O\left(\left(\hat\lambda / \lambda_1\right)^{i} + \left(\hat\nu / \nu_1\right)^{j}\right)$.
Now we will select  $\mu_2 = O\left(\left(\hat\lambda / \lambda_1\right)^{i} + \left(\hat\nu / \nu_1\right)^{j}\right)$ in the way to make 
the value of $\delta_{ij} = \mu_2 + \rho_{ij}^3$ asymptotically small as $i,j \to \infty$. Then we have 
$ A_{21}^{(2)}(\mu_2 + \rho_{ij}^5) =  A_{21}^{(2)} \rho_{ij}^6 = 
O\left(\left(\hat\lambda / \lambda_1\right)^{i} + \left(\hat\nu / \nu_1\right)^{j}\right) \neq 0$
as otherwise the Jacobian of map $T_{21} T_{02}^jT_{12} T_{01}^i$ would be vanishing when $\delta_{ij}$ goes to zero.

Finally  we rescale the coordinates as follows
$$
\displaystyle  v = -\frac{\gamma_1^{-i}\gamma_2^{-j}}{d^{(1)} d^{(2)}} (1 + q_{ij}^{(3)}) \; Y \;,\; 
x_1  = - \frac{b_1^{(2)} \gamma_1^{-i} \gamma_2^{-j}} {d^{(1)} d^{(2)}}\;X_1 \;,\; 
x_2  = - \frac{b_1^{(2)} A_{21}^{(2)}\rho_{ij}^6 \lambda_1^i \gamma_1^{-i} \lambda_2^j \gamma_2^{-j}}{d^{(1)} d^{(2)}}\;X_2 \;.
$$
After this, we can rewrite (\ref{tk-1+n}) in the following form
\begin{equation}
\label{tkres31}
\begin{array}{l}
 \bar X_1 \; = \;  Y + O( \lambda_1^i \nu_1^j  + \gamma_1^{-i} \gamma_2^{-j}))\;, \\
\bar X_2  \; = \;  X_{1} + O( \lambda_1^i \nu_1^j  + \gamma_1^{-i} \gamma_2^{-j})) \;,\\
 \bar Y = M_1 + M_2 X_1 + J_{ij} X_2 - Y^2 + O( \lambda_1^i \nu_1^j  + \gamma_1^{-i} \gamma_2^{-j})).
\\
\end{array}
\end{equation}


\section*{ACKNOWLEDGMENTS}
the author thanks J.\,C. Tatjer for useful discussions that led to the idea of condition {\bf D2} and {\bf Case~{\rm III}}.

\section*{FUNDING}
This paper is a contribution to the project M7 (Dynamics of Geophysical Problems in Turbulent Regimes) of the Collaborative Research Centre TRR 181 ``Energy Transfer in Atmosphere and Ocean'' funded by the Deutsche Forschungsgemeinschaft (DFG, German Research Foundation) -- Projektnummer 274762653. The paper is also supported by the
grant of the Russian Science
Foundation 19-11-00280.

\section*{CONFLICT OF INTEREST}
The author declares that there are no conflicts of interest.

\section*{SUPPLEMENTARY MATERIALS}
The data that support the findings of this study (proofs, pictures) are placed in the body of the
text. If some extra requirements appear, they should be addressed to the corresponding author.

\endpaper

\end{document}